\magnification 1200
\input amssym.def
\input amssym.tex
\parindent = 40 pt
\parskip = 12 pt

\font \heading = cmbx10 at 12 true pt
 at 22 true pt
\font \medheading =cmbx7 at 16 true pt
 at 7 true pt

\def \R{{\bf R}}

\centerline{\medheading Convolution kernels of 2D Fourier  multipliers}
\centerline{\medheading based on real analytic functions}
\rm
\line{}
\line{}
\centerline{\heading Michael Greenblatt}
\line{}
\centerline{June 26, 2016}
\baselineskip = 12 pt
\font \heading = cmbx10 at 15 true pt
 at 13 true pt
\line{}
\line{}
\line{} 
\noindent{\heading 1. Background and Theorem Statements}.

\vfootnote{}{This research was supported in part by NSF grant  DMS-1001070.}

In this paper, estimates are proven for convolution kernels associated to multipliers
from a reasonably general class of compactly supported two-dimensional functions constructed 
out of real analytic functions. These estimates are both for overall decay rate and decay 
rate in specific directions. The estimates are sharp for a certain range of exponents 
appearing in the theorems. In a separate paper [G3], a class of 
"well-behaved" functions is described that contains a number of relevant examples and for which,
after a little more work, 
these estimates can be explicitly described in terms of the Newton polygon of the function.

The  compactly supported  Fourier multipliers $m(x,y)$ we consider are as follows. For each 
$(x_0,y_0)$ in the support of $m(x,y)$ we assume that on a neighborhood
of $(0,0)$  the function $m(x_0 + x, y_0 + y)$ can be written in the form
$$m(x_0 + x, y_0 + y) = \alpha(x,y)\, \chi_E(x,y) \,\prod_{i = 1}^n |f_i(x,y)|^{\gamma_i}  \eqno (1.1)$$
Here $\alpha(x,y)$ is $C^1$ except at $(0,0)$ and for some constant $A$ one has
$$|\alpha(x,y)| \leq A \eqno (1.2)$$
$$|\nabla \alpha(x,y)| \leq A (x^2 + y^2)^{-{1 \over 2}} \eqno (1.3)$$
The prototypical $\alpha(x,y)$ would be a smooth function supported on a neighborhood of the origin.
The functions $f_i(x,y)$ are real analytic and not identically zero on a neighborhood of the origin. The set $E$ is assumed to be either a disk $\{(x,y): x^2 + y^2 < r^2\}$ or to be expressible as a disjoint union  of open sets
$\cup_{i = 1}^m E_i$, where each $E_i$ is a region bounded by curves $C_1$, $C_2$ connecting the origin to a circle
 $x^2 + y^2  = r^2$, and the circle $x^2 + y^2  = r^2$ itself. The curves
$C_1$ and $C_2$  are assumed to be either half of the graph of the form $y = h(|x|^{1 \over N})$ or $x = h(|y|^{1 \over N})$ for a real analytic
$h$ with $h(0) = 0$. There are two regions formed by the curves $C_1$, $C_2$, and the circle and we allow a given $E_i$ to be either of them. We assume that all the curves
$C_1$ and $C_2$ are disjoint in the disk $\{(x,y): x^2 + y^2 \leq r^2\}$ so that the $E_i$ are wedge or sliver-shaped regions
whose closures only intersect at the origin.

The above form of $E$ is a convenient way to describe a general domain defined through real analytic functions. In fact 
any such curve is part of the zero-set of a real analytic function. For example, in the case of the graph of $y = h(x^{1 \over N})$
one can take $\prod_{j=0}^{N-1}(y -  h(e^{2\pi ij \over N}x^{1 \over N}))$. Conversely, by Puiseux's theorem the zero
set of a real analytic function is locally the finite union of curves of the form used here.

The form of $E$ used here allows us for example  to define the multiplier in several ways on several regions. If the different regions can be defined via real analytic functions, then one can write the multiplier as the sum of several multipliers of the form used here, and then
add the kernel estimates obtained by our theorems. 
Another reason to use this form is if instead of  wanting $|f_i(x,y)|^{\gamma_i}$ in the multiplier, you wanted a factor to reflect the sign of $f_i(x,y)$, then you could write the multiplier as the 
sum of two terms depending on the sign of $f_i(x,y)$; the curves where $f_i(x,y) = 0$ can be incorporated into the boundary of 
$E$.

The only restriction we assume on the exponents $\gamma_i$ is that $\chi_E(x,y) \prod_{i=1}^n f_i(x,y)^{\gamma_i}$ is integrable
on a neighborhood of the origin; otherwise even taking the Fourier transform of $m(x,y)$ would involve delicate distribution theory
issues.

Using a partition of unity we can write $m(x,y) = \sum_{j=1}^K m_j(x, y)$,
where each $m_j(x_j + x, y_j + y)$ satisfies $(1.1)$ for some $(x_j,y_j)$. The convolution kernel of $m(x,y)$ can then be written in the form
$$K( t,u) = \sum_{j=1}^M \int  m_j(x, y) e^{it x + iu y}\,dx\,dy \eqno (1.4)$$
Although the multipliers of this paper do not appear to have been extensively studied before, they are related to damped scalar
oscillatory integrals of the form
$$G(s,t,u) = \int_{\R^2} |f(x,y)|^{\alpha}e^{is S(x,y) + itx + iuy} \phi(x,y)\,dx\,dy $$
Here $\phi(x,y)$ is a smooth cutoff function supported near the origin, $S(x,y)$ is a real analytic function near the origin with $S(0,0) = 0$ and $\nabla S(0,0) = (0,0)$, and one seeks estimates of 
the form $|G(s,t,u)| \leq C(1 + |(s,t,u)|)^{-\epsilon}$. By taking $s = 0$ one is reduced to situations studied in this paper.
Such oscillatory integrals come up frequenly when using the damping function techniques initiatied in [SoS] when studying maximal
 averages over surfaces, such as in the papers [CMa1] [IM] [IoSa1] [IoSa2] [G4]. On their own, such 
oscillatory integrals can be viewed as surface measure Fourier transforms for surfaces with damping functions, possibly singular,
such as those considered in [CDMaM] [CMa2] [G1] [Gr].

In order to state the main theorems of this paper, we will need a couple of facts following from resolution of singularities which we
will prove at the end of section 2.

\noindent {\bf Lemma 1.1.} Let $g(x,y) = \chi_E(x,y) \,\prod_{i = 1}^n |f_i(x,y)|^{\gamma_i}$, where $E, f_i,$ and 
$\gamma_i$ are as before.  There exist $c_1,c_2,c_3 > 0$, an $\epsilon > 0$, and $d = 0$ or $1$ such that if $0 < r < c_3$ one has
$$c_1 r^{\epsilon} |\ln r|^d \leq \int_{x^2 + y^2 < r^2}  g(x,y)\,dx\,dy \leq  c_2 r^{\epsilon} |\ln r|^d \eqno (1.5)$$
\noindent {\bf Lemma 1.2.}  Let $g(x,y)$ be as in Lemma 1.1. Suppose $v = (v_1,v_2)$ is a unit vector in $\R^2$, and let 
$v^{\perp} = (v_2, -v_1)$ be the  orthogonal unit vector.  There exist  $\delta_v, c_v > 0$ and a $e_v = 0$ or
 $1$ such that if $c < c_v$ then there are  $a_{v,c},b_{v,c} > 0$ such that for $0 < r < c$ one has 
$$a_{v,c} r^{\delta_v}|\ln r|^{e_v} \leq \int_{\{(x,y): |(x,y) \cdot v| < r,\, |(x,y) \cdot v^{\perp}| < c\}} g(x,y)\,dx\,dy \leq b_{v,c} r^{\delta_v}|\ln r|^{e_v} \eqno (1.6)$$
Note that for any direction $v$, the rate of decrease in $(1.5)$ is at least as fast as the decrease rate in $(1.6)$ since the 
domain of integration in $(1.6)$ contains the disk of radius $r$ centered at the origin, which is the domain of integration in $(1.5)$.

We now give the local theorems for the (inverse) Fourier transform of $m(x,y)$ which will sum to give
 the overall kernel estimates.
We use the following notation. Let $\phi(x,y)$ be a nonnegative bump function which is one on a neighborhood of $(x_0,y_0)$, and
let $m_{\phi, x_0, y_0}(x,y) = \phi(x_0 + x,y_0 + y) m(x_0 + x, y_0 + y)$. We assume that the support of $\phi(x,y)$ is small enough so that $\phi(x_0 + x,y_0 + y) m(x_0 + x, y_0 + y)$ can be written in the form $(1.1)$. We then define $K_{\phi,x_0,y_0}(t,u)$ by
 $$ K_{\phi,x_0,y_0}(t,u) = \int_{\R^2}\phi(x,y) m(x,y)e^{itx + iuy}\,dx\,dy$$
$$= e^{itx_0 + iuy_0}\int_{\R^2} m_{\phi, x_0, y_0}(x,y)e^{itx + iuy}\,dx\,dy$$
 Thus $K_{\phi,x_0,y_0}(t,u)$ can be viewed as the contribution to the convolution kernel of 
$m(x,y)$ coming from the region near $(x_0,y_0)$.

For each $f_i(x,y)$ appearing in $(1.1)$,
let $F_i(x,y)$ be the sum of the terms of $f_i(x,y)$'s Taylor expansion at $(0,0)$ of lowest total degree. The zeroes of a given 
$F_i(x,y)$ are either a finite union of lines through the origin, just the origin, or the empty set (in the case when $f_i(0,0) \neq 0$). We
let $l_1,...l_{p'}$ be the list of all such lines over all $i$ (if there are any). We add to this list any lines that are tangent at the
origin to the boundary curves $C_1$ and $C_2$ of the $E_i$ as described after $(1.3)$. We denote the combined list of lines by $l_1,...,l_p$, with 
the understanding that the combined list might be empty.

\noindent We get the strongest results when the $\epsilon$ in Lemma 1.1 is less than ${1 \over 2}$:

\noindent {\bf Theorem 1.3.} Suppose $(x_0,y_0)$ is in the support of $m(x,y)$ and let $\epsilon$ and $d$ be as in Lemma 1.1 as applied to the $g(x,y)$ associated with
 $m(x_0 + x, y_0 + y)$. If $\epsilon < {1 \over 2}$, then the following hold, where
$|(t,u)|$ denotes the magnitude $(t^2 + u^2)^{1 \over 2}$ of the vector $(t,u)$.

\noindent {\bf a)} For a given line $l$ through the origin, let $l_H$ denote the points in $\R^2$ within distance $H$ of $l$. Let
$\delta_v$ and $e_v$ be as in Lemma 1.2, where $v$ is in the direction of $l$. If $v$ is perpendicular to one of the lines
$l_1,...,l_p$, then if the support of $\phi(x,y)$ is sufficiently small, depending on $v$, 
there is a constant $C$ depending $g(x,y)$, $\phi(x,y)$,
$H$, $l$, and the constant $A$ of $(1.2)-(1.3)$ such that for $(t,u)$ in the strip  $l_H$ one has 
$$|K_{\phi,x_0,y_0}(t,u)| \leq C(2 + |(t,u)|)^{-\delta_v }(\ln (2 + |(t,u)|))^{e_v} \eqno (1.7a)$$
If $v$ is not perpendicular to a $l_i$, then $\delta_v = \epsilon$ and instead of $(1.7a)$ we have the estimate
$$|K_{\phi,x_0,y_0}(t,u)| \leq C(2 + |(t,u)|)^{-\epsilon }(\ln (2 + |(t,u)|))^d \eqno (1.7b)$$

\noindent {\bf b)} Let $(\delta, e)$ denote  the slowest decay rate in part a) over all lines.
If the support of $\phi(x,y)$ is
sufficiently small there is a constant $C'$ depending on $g(x,y)$, $\phi(x,y)$, and $A$ such that for any $(t,u)$ one  has the estimate
$$|K_{\phi,x_0,y_0}(t,u)| \leq C'(2 + |(t,u)|)^{-\delta }(\ln (2 + |(t,u)|))^e \eqno (1.8)$$

\noindent {\bf c)} If there exists a $c > 0$ such that  $\alpha(x,y)$ in $(1.1)$ satisfies $\alpha(x,y) > c$ on a  neighborhood of the origin, then
parts a) and b) of this theorem are sharp in the sense that the exponents $\delta_v$, $\delta_v = \epsilon$, and
$\delta$ cannot be improved in $(1.7a)$, $(1.7b)$, and $(1.8)$ respectively.

When $\epsilon \geq {1 \over 2}$ but some $\delta_v < {1 \over 2}$,  we have the following weaker 
version of Theorem 1.3, which still gives the
 optimal overall decay rate of part b), but which does not give the best estimates in all directions.

\noindent {\bf Theorem 1.4}. Suppose $(x_0,y_0)$ is in the support of $m(x,y)$ and let $\epsilon$ and $d$ be as in Lemma 1.1 as applied to the $g(x,y)$ associated with
 $m(x_0 + x, y_0 + y)$. If $\epsilon \geq {1 \over 2}$, but there is at least one direction for which $\delta_v < {1 \over 2}$,
 then the following hold.

\noindent {\bf a)} There are at most finitely many directions for which the corresponding  $\delta_v$ is less than ${1 \over 2}$, and 
each such direction  must be perpendicular to one of the lines $l_1,...,l_p$. For each such direction,
we have the same estimate as in Theorem 1.3:  if the support of $\phi(x,y)$ is sufficiently small there is a constant $C$ depending $g(x,y)$, $\phi(x,y)$,
$H$, $l$, and $A$ such that for $(t,u)$ in the strip $l_H$ one has
$$|K_{\phi,x_0,y_0}(t,u)| \leq C(2 + |(t,u)|)^{-\delta_v}(\ln (2 + |(t,u)|))^{e_v} \eqno (1.9)$$
This estimate is sharp in the same sense as in Theorem 1.3 c).

\noindent For the remaining directions, we still have the (usually nonsharp) estimate that in place of $(1.9)$ one has
$$|K_{\phi,x_0,y_0}(t,u)| \leq C(2 + |(t,u)|)^{-{1 \over 2}}(\ln (2 + |(t,u)|))^2 \eqno (1.10)$$
\noindent {\bf b)} The statement of part b) of Theorem 1.3 holds and is sharp in the same sense as in Theorem 1.3.

Our next theorem says that in the case that all $\delta_v$ are at least ${1 \over 2}$, one still gets an exponent of at least ${1 \over 2}$ in any direction, and also for
the overall decay rate. As a result, Theorems 1.3 and 1.4 give the best overall decay rate whenever it is less than
${1 \over 2}$. 

\noindent {\bf Theorem 1.5.} Let $\epsilon$ and $d$ be as in Theorems 1.3 and 1.4.
 If $\delta_v \geq {1 \over 2}$ for all directions $v$, then there is a constant
 $C$ depending on $g(x,y)$, $\phi$, and $A$ such that one has the estimate
$$|K_{\phi,x_0,y_0}(t,u)| \leq C(2 + |(t,u)|)^{-{1 \over 2}}(\ln (2 + |(t,u)|))^2 \eqno (1.11)$$
The above theorems give local estimates for the convolution kernel associated to a given $m(x,y)$ of the type treated in the paper. 
One can then use a partition of unity to write $m(x,y) = \sum_{i=1}^K m_i(x,y)$, where one of the above theorems provides
estimates for each $m_i(x,y)$, thereby giving global estimates for this kernel. When one obtains a sharp estimate for any $m_i(x,y)$,
one typically obtains a sharp estimate for $m(x,y)$ as well; cancellation does not typically occur. We describe this phenomenon in the
next theorem.

\noindent {\bf Theorem 1.6.} Suppose $m(x,y) = \sum_{i=1}^K m_i(x,y)$ such that each $m_i(x,y)$ is localized enough so that
Theorem 1.3 or Theorem 1.4 applies to $m_i(x,y)$. Suppose there is a $c > 0$ such that each function 
$\alpha(x,y)$ of $(1.1)$ corresponding to any $m_i(x,y)$ satisfies $\alpha(x,y) > c$ on the support of $m_i(x,y)$. Suppose further
that when
adding the estimates given by Theorems 1.3 or 1.4 the resulting estimate is one that is stated by Theorem 1.3 or 1.4 to be
 sharp for at least one of the $m_i(x,y)$ that it came from. Then this
estimate is also sharp for $m(x,y)$ in the same sense that it was stated to be sharp for any such $m_i(x,y)$.  

To help understand heuristically why in general one will not get a better exponent than ${1 \over 2}$ than in Theorems 1.3-1.6, 
we focus on Theorems 1.3a) and 1.4a) and 
 consider the case where $E = \{(x,y) \in D: x > 0,\, x^2 < y < 2x^2\}$, 
where $D$ is a small disk centered at the origin, and assume there are two $f_i(x,y)$, given by $f_1(x,y) = x$ and $f_2(x,y) = y - x^2$.
We make no restrictions on $\gamma_1$, and let $\gamma_2 = -1 + \eta$ for some small $\eta$.
 Assume $\alpha(x,y)$ is identically equal to $1$.
Then the convolution kernel associated to the multiplier in this case is given by
$$K(t,u) = \int_D x^{\gamma_1}(y - x^2)^{-1 + \eta}e^{itx + iuy}\,dx\,dy \eqno (1.12)$$
Changing variables from $y$ to $y + x^2$ and setting $t = 0$, we get
$$K(0,u) = \int_{\{(x,y) \in D:\, x > 0, \,\,0 < y < x^2\}}x^{\gamma_1} y^{-1 + \eta}e^{iux^2 + iuy} \eqno (1.13)$$
When $\eta$ is very small, the $y^{-1 + \eta}$ factor ensures that one gets very little decay in $K(t,u)$ due to the $iuy$ term in 
the exponential; the behavior is driven by the $x$ integral in $(1.13)$ for fixed values of $y$. Stationary phase can be readily 
used on each dyadic piece of this $x$ integral and the result is
$$|K(0,u)| \leq C \int_{\{(x,y) \in D:\, x > 0, \,\,0 < y < x^2\}}x^{\gamma_1} y^{-1 + \eta}\min\bigg(1, 
{1 \over |ux^2|^{1 \over 2}}\bigg) \eqno (1.14)$$
Converting back to the original variables and using that $y \sim x^2$ on the domain of integration yields
$$|K(0,u)| \leq C \int_{\{(x,y) \in D:\, x > 0, \,\,x^2 < y < 2x^2\}} x^{\gamma_1}(y - x^2)^{-1 + \eta}\min\bigg(1, 
{1 \over |uy|^{1 \over 2}}\bigg) \eqno (1.15)$$
Because of the exponent ${1 \over 2}$ in the ${1 \over |uy|^{1 \over 2}}$ factor in $(1.15)$, in the $u$ direction one can never
 get a better decay 
rate than $|u|^{-{1 \over 2}}$ in $(1.15)$. The $u$ direction here corresponds to a direction perpendicular to a $l_i$ in 
Theorems 1.3-1.4. At the same time, one may select $\gamma_1$ such that the exponent $\delta_v$ in Theorem 1.3-1.4  is a given
value greater than ${1 \over 2}$. While there is a slight improvement over the above heuristics due to the  $iuy$ term in $(1.13)$, as $\eta$ goes to zero,
this improvement vanishes.  Hence the statements of Theorems 1.3a) and 1.4a) will not hold in
 generality if we replace  ${1 \over 2}$ by any larger exponent. 
Similar considerations apply concerning the optimality of this exponent in the other parts of Theorems 1.3-1.6.

Examples like the above show that the sharp estimates of Theorems 1.3 and 1.4 do not hold in general if the exponents are greater 
than ${1 \over 2}$. However, the sharpness proofs we will give in section 4 do extend to any $\delta_v > 0$ and $\epsilon > 0$ situations, meaning that in such situations one cannot prove better estimates than the above sharp estimates either.  It is unclear if there is
a general statement that can be stated simply that covers index ranges beyond ${1 \over 2}$. We will however prove a 
theorem which does give at least some estimates in these ranges:

\noindent {\bf Theorem 1.7.} Suppose $1 < p \leq \infty$ is such that for each $(x_0,y_0)$ in the support of $m(x,y)$, the function 
$g(x,y) =  \chi_E(x,y) \,\prod_{i = 1}^n |f_i(x,y)|^{\gamma_i}$ is in $L^p(N)$ for some neighborhood $N$ of the origin. Then if
$p'$ denotes the complementary exponent satisfying ${1 \over p} + {1 \over p'} = 1$, for some constant $C$ depending on $m(x,y)$ and $p$ one has that $|K(t,u)| \leq C(2 + |(t,u)|)^{-{1 \over p'}}$ when $p < \infty$, and $|K(t,u)| \leq C(2 + |(t,u)|)^{-1}
\ln(2 + |(t,u)|)$ if $p = \infty$.

\noindent {\heading 2. Resolution of singularities in two dimensions and some consequences.}

We will make use of the real-analytic case of the resolution of singularities theorem of [G5], which is an extension of related theorems in [G1]-[G2], and which was also influenced by [PS] and [V]. This theorem is as follows. Let $S(x,y) = 
\sum_{\alpha,\beta} S_{\alpha\beta}x^{\alpha}y^{\beta}$ be a real-analytic
function on a neighborhood of the origin, not identically zero, satisfying $S(0,0) = 0$.

Divide the
$xy$ plane into eight triangles by slicing the plane using  the $x$ and $y$ axes  and two lines through the origin, one of the form $y = mx$ for some $m > 0$ and one of the form $y = mx$ for some $m < 0$. One must ensure that these two lines are not ones  on 
which the function $\sum_{\alpha + \beta = o} S_{\alpha\beta}x^{\alpha}y^{\beta}$ vanishes other than at the
origin, where $o$ denotes the order of the zero of $S(x,y)$ at the origin. After reflecting about the $x$ and/or $y$ axes and/or the line $y = x$ if necessary, each of the triangles becomes of the form $T_b = \{(x,y) \in \R^2: x > 0,\,0 < y < bx\}$ (modulo an inconsequential boundary set of measure zero). We first give a 
relevant version of Theorem 2.1 of [G1].

\noindent {\bf Theorem 2.1.} Let  $T_b = \{(x,y) \in \R^2: x > 0,\,0 < y < bx\}$ be as above. Abusing notation slightly, use the notation $S(x,y)$ to denote the reflected function $S(\pm x,\pm y)$ or $S(\pm y, \pm x)$ corresponding to $T_b$.
 Then there is a $a > 0$ and a positive integer $N$ such that
if $F_a$ denotes  $\{(x,y) \in \R^2: 0 \leq x\leq a, \,0 \leq y \leq bx\}$, then one can write $F_a = \cup_{i=1}^n cl(D_i)$, such that for to each $i$ there is a $k_i(x) = l_i x^{s_i} + ...$ with $k_i(x^N)$ real-analytic and $s_i \geq 1$ such that after a coordinate change of the form $\eta_i(x,y) = (x, \pm y + k_i(x))$, the set $D_i$ becomes a set $D_i'$ on which the function $S \circ \eta_i(x,y)$ approximately becomes a monomial $d_i x^{\alpha_i}y^{\beta_i}$, $\alpha_i$ a nonnegative rational number and $\beta_i$ a nonnegative integer in the following sense.

\noindent {\bf a)} $D_i' = \{(x,y): 0 < x < a, \, g_i(x) < y < G_i(x)\}$, where $g_i(x^N)$ and $G_i(x^N)$ are
real-analytic. If we expand $G_i(x) =  H_i x^{M_i} + ...$, then $M_i \geq 1$ and $H_i > 0$, and consists of a single term $H_ix^{M_i}$ when $\beta_i = 0$. 

\noindent {\bf b)} Suppose $\beta_i = 0$. Then $g_i(x) = 0$. The set $D_i'$ can
be constructed such that for any predetermined $\eta > 0$ there is a $d_i \neq 0$ such that on $D_i'$, for all $0 \leq l \leq \alpha_i$ one has
$$ |\partial_x^l (S \circ \eta_i)(x,y) -  d_i\alpha_i(\alpha_i -1) ... (\alpha_i - l + 1)x^{\alpha_i - l}| < \eta|d_i| x^{\alpha_i-l} \eqno (2.1)$$ 

\noindent {\bf c)} If $\beta_i > 0$, then $g_i(x)$ is either identically zero or $g_i(x)$ 
can be expanded as $h_ix^{m_i} + ...$ where $h_i > 0$ and $m_i > M_i$. The $D_i'$ can
be constructed such that such that for any predetermined $\eta > 0$ there is a $d_i \neq 0$ such that on $D_i'$, for all $0 \leq l \leq \alpha_i$ and all $0 \leq  m \leq \beta_i$ one has
$$|\partial_x^l\partial_y^m(S \circ \eta_i)(x,y) -  \alpha_i(\alpha_i - 1)...(\alpha_i - l + 1)\beta_i(\beta_i - 1)...(\beta_i - m + 1)
d_ix^{\alpha_i - l}y^{\beta_i - m}| $$
$$\leq \eta |d_i| x^{\alpha_i - l}y^{\beta_i - m} \eqno (2.2)$$
It should be pointed out that in Theorem 2.1 of [G1] (but not Theorem 3.1 of [G2]) it was assumed that one had rotated coordinates
in advance so that $\partial_x^o S(0,0) \neq 0$ and $\partial_y^o S(0,0) \neq 0$, where $o$ is the order of the zero of $S(x,y)$ at
the origin. This was done to make the exposition of the smooth situation somewhat easier, and is not necessary for the arguments to work.

For the purposes of proving our theorems, we will need to simultaneously resolve the singularities of several functions.
As is well-known in the 
subject of resolution of singularities, one can often simultaneously resolve the singularities of several functions by resolving the 
singularities of their product. This was done in [G5], where the following theorem was proven.

\noindent {\bf Theorem 2.2.} Suppose $S_1(x,y),...,S_k(x,y)$ are real-analytic functions on a neighborhood of the origin, none
identically zero, with
$S_j(0,0) = 0$ for each $j$. Let $D_i'$, $\alpha_i$, and $\beta_i$ be as in Theorem 2.1 applied to  $\prod_{j=1}^k S_j(x,y)$.
 Then one can further divide each $D_i'$ into finitely many pieces $D_{il}'$, such that on each $D_{il}'$ an
additional coordinate change of the form $(x,y) \rightarrow (x,y - c_{il} x^{M_i})$ or $(x,y -  c_{il}x^{m_i})$, $c_{il} \geq 0$, will
result in each $S_j(x,y)$ satisfying the conclusions of Theorem 2.1, with one difference: the domains $D_{il}'$ with $\beta_i = 0$
now are only assumed to have the same form as the domains where $\beta_i > 0$. That is, $D_{il}'$ is the form $\{(x,y): 0 < x < a, \, g_{il}(x) < y < G_{il}(x)\}$, where $g_{il}(x^N)$ and $G_{il}(x^N)$ are real-analytic,  $G_{il}(x) =  H_{il} x^{M_{il}} + ...$, and  $g_{il}(x) = h_{il} x^{m_{il}} + ... $ where $1 \leq M_{il} < m_{il}$ and $h_{il} \geq 0,  H_{il} > 0$.

\noindent We will also use the following corollary to Theorem 2.2 which was proven in [G5].

\noindent {\bf Corollary 2.3.} For any given $K$, however large, for any predetermined $\eta > 0$ the $D_i'$ can be constructed so that $(2.1)$ and $(2.2)$  hold for all $\alpha_i, \beta_i < K$.

\noindent {\bf Proof of Lemma 1.1.}  

If each $f_i(0,0) \neq 0$, the result easily follows by finding
the area of the portion of $E$ within distance $r$ of the origin, so we assume at least one $f_i(0,0) \neq 0$. We can also replace 
each $f_i$ for which $f_i(0,0) \neq 0$ by the constant function $1$, so without loss of generality we can remove these functions and
assume that each $f_i(0,0) = 0$. 

It suffices to prove $(1.5)$ replacing integrals over
discs centered at the origin with integrals over rectangles of fixed edge length ratio, and this is what we will do. We apply
Theorem 2.2 to $f_1,...,f_n$, and the result is a rectangle centered at the orgin on which Theorem 2.2  holds. We will show $(1.5)$ 
for dilations of this rectangle. Theorem 2.2 provides slivers of the form 
$S = \{(x,y): 0 < x < a,  g_{il}(x) < y < G_{il}(x)\}$ with $g_{il}(x^N)$ and  $G_{il}(x^N)$ real analytic for some
positive integer $N$. On this set, in the new coordinates each $|f_i(x,y)|$ is within a constant factor of 
some $x^{\alpha_i}y^{\beta_i}$. Thus the
 product $\prod_{i=1}^n |f_i(x,y)|^{\gamma_i}$ is also within a constant factor of some $M(x,y)= 
x^{\alpha}y^{\beta}$.

If one integrates $M(x,y)$ over the set $\{(x,y): 0 < x < r,\,  g_{il}(x) < y < G_{il}(x)\}$, one obtains an expression of the form $c r^a(\ln r)^b + o(r^a(\ln r)^b)$, where the $r^a(\ln r)^b$ term is derived from the
leading terms of the Taylor expansions of $ g_{il}(x)$ and $G_{il}(x)$. Here $b = 0$ or $1$. Since the coordinate changes of Theorem 2.2 all have Jacobian
1, the integral of $M(x,y)$ over this sliver in its original coordinates will be of the same form.

 If one
now inserts a $\chi_E(x,y)$ factor and  looks at the integral of $\chi_E(x,y) M(x,y)$ over the sliver $S$ in the original coordinates, and 
transfers to the new coordinates, instead of integrating over $S$ in the new coordinates, one integrates over a portion 
 cut out by at most
finitely many functions of the form $y = g(x)$ where some $g(x^N)$ is real analytic. Again direct integration reveals that
the result is of the same form $c r^a(\ln r)^b + o(r^a(\ln r)^b)$. Hence the integral of $M(x,y)$ over $S \cap  E$, in the
original coordinates or final coordinates, is of this form. Since $\chi_E(x,y) \prod_{i=1}^n |f_i(x,y)|^{\gamma_i}$ is within a bounded
factor of $M(x,y)$, we conclude that  the integral of $\chi_E(x,y) \prod_{i=1}^n |f_i(x,y)|^{\gamma_i}$
is also within a constant factor of some $c r^a(\ln r)^b$. Adding this over all slivers gives $(1.5)$, completing the proof of Lemma 1.1.

\noindent {\bf Proof of Lemma 1.2.}

Let $l$ be a line segment centered at the origin with direction $v^{\perp}$ such that each
$f_i(x,y)$ and each of the functions defining $E$ is defined on a neighborhood of $l$. 
Using a partition of unity, we let $p_1,...,p_k$ be points on $l$ such that to each $p_j$ there is a rectangle $R_j$ 
centered at $p_j$ such that either the product $\prod_i^n f_i(x,y)$ is nonzero on a neighborhood of $cl(R_j)$ or such that 
Theorem 2.2 holds for the product of the nonzero $f_i(x,y)$ on the rectangle $R_j$ when we center at $p_j$ and have rotated so 
that the $v^{\perp}$ direction has become the $x$ direction. It suffices to prove $(1.6)$ for 
the portion of the integral contained in a given $R_j$ since the overall result will follow simply by adding these statements over all $j$.

As in part a), the estimates for the rectangles where $\prod_i^n f_i(x,y)$ is nonzero follow from a straightforward integration, so
we assume at least one $f_i(x,y)$ is zero at $p_j$. Analogous to  part a) we can assume the partition of unity is such that we
may replace all of the $f_i(x,y)$ which are nonzero at $p_j$ by the constant function 1. Thus without loss of generality we can 
assume each $f_i(p_j) = 0$. Since we have rotated so that $v^{\perp}$ is the $x$ direction, our goal is to understand as a function of
$r$ the integral
of $\chi_E(x,y) \prod_{i=1}^n |f_i(x,y)|^{\gamma_i}$ over the portion of $R_j$ for which $|y| < r$. As in part a), it suffices 
to show $(1.6)$ for the portion of the integral over $|y| < r$ coming from each of the slivers arising from Theorem 2.2, as the
overall result will then follow via adding over all slivers.

If the sliver is one of the ones adjacent to the upper or lower boundaries of the rectangle $R_j$, then the coordinate changes 
of Theorem 2.2 turn the lines $y 
= \pm r$ into the line $x = r$, and the situation reduces to the one considered in part a), so we have the desired estimates in this
situation. Assume therefore that
the sliver is one of the ones adjacent to the right or left boundaries of $R_j$. The overall coordinate change in Theorem 2.2 
is of the form $(x, y) \rightarrow (\pm x, \pm y + k(x))$, where some $k(x^N)$ is real analytic. If $k(x)$ happens to be the
zero function, then $\prod_{i=1}^n |f_i(x,y)|^{\gamma_i}$ is already comparable in magnitude to some $M(x,y)$ of the
form $x^cy^d$, so
one may perform a direct integration of $M(x,y)$ to get an expression of the form $c r^a(\ln r)^b + o(r^a(\ln r)^b)$.
The presence of a $\chi_E(x,y)$ factor will not change the resulting form, for the same reasons as in part a).

If $k(x)$ is not identically zero, we denote by $p$ the degree of the initial term of the Taylor expansion of $k(x)$ at the origin. Cutting 
off the sliver at height $y = r$ or $y = -r$ in the original coordinates has a similar  effect as cutting off the sliver with a vertical line 
$x = r^{1 \over p}$ or $- r^{1 \over p}$; when $k(x)$ is not identically zero, by construction 
the sliver in the original coordinates is always contained within a wedge $c_1 |x|^p < y < c_2 |x|^p$ that is in one of the
four quadrants. In view of the monomial form of the functions in the final coordinates,  the integral of $\chi_E(x,y) \prod_{i=1}^n |f_i(x,y)|^{\gamma_i}$ over the portion of the sliver in the original coordinates where
$y < r$ will therefore be within a constant factor of the integral over the portion of the sliver in the original coordinates where $x < r^{1 \over p}$. This can then
be computed directly in the same way as one computed the integral over the portion where $x < r$ for the first kind of 
sliver, and in part a) of this
lemma. The result will once again be comparable to $ r^a(\ln r)^b$ for some $a$ and $b$. Thus we see that we have such
 an expression for all slivers, and the proof of Lemma 1.2 is complete.

\noindent {\heading 3. Proofs of the estimates of Theorem 1.3, 1.4, 1.5 and Theorem 1.7.}

\noindent We start with the well-known Van der Corput lemma (see p. 334 of [S]).

\noindent {\bf Lemma 3.1.} Suppose $k \geq 2$ and $h(x)$ is a $C^k$ function on the interval $[a,b]$ with $|h^{(k)}(x)| > A$ on $[a,b]$ for
some $A > 0$. Let $\phi(x)$ be $C^1$ on $[a,b]$. If $k \geq 2$ there is a constant $c_k$ depending only on $k$ such that
$$\bigg|\int_a^b e^{ih(x)}\phi(x)\,dx\bigg| \leq c_kA^{-{1 \over k}}\bigg(|\phi(b)| + \int_a^b |\phi'(x)|\,dx\bigg) $$
If $k =1$, the same is true if we also assume that $h'(x)$ is monotonic on $[a,b]$. 

Throughout most of this section, we will be focusing on local behavior near a given $(x_0,y_0)$. Namely, using the notation of
$(1.1)$, for 
$\phi(x,y)$ supported on a small neighborhood of $(x_0,y_0)$ and various sets $S$ we will be looking at quantities of the form
$$\bigg|\int_S \phi(x_0 + x,y_0 + y) \alpha(x,y)\, \chi_E(x,y) \,\prod_{i = 1}^n |f_i(x,y)|^{\gamma_i} e^{itx + iuy}\,dx\,dy\bigg|$$
To simplify notation, we will just write $\alpha(x,y)$ in place of  $\phi(x_0 + x,y_0 + y)\alpha(x,y)$ with the understanding that $\alpha(x,y)$
is to be supported on a sufficiently small neighborhood of $(0,0)$ for our arguments to work.

Our next lemma provides the key Fourier transform estimate for a given sliver arising from Theorem 2.1 or Theorem 2.2. Theorems 
1.3 and 1.4 will
be proven by adding these estimates over all slivers and interpreting the result in an appropriate way.

\noindent {\bf Lemma 3.2.} Let $S$ be a sliver in the original coordinates arising from an application of Theorem 2.2  to
$f_1(x,y)$,...,$f_n(x,y)$, and real analytic functions whose zero sets contain all the boundary curves of $E$ on a 
neighborhood of the origin (recall such functions always exist).
  Then if the function $\alpha(x,y)$ in $(1.1)$ is supported on the
 neighborhood of the origin on which we are applying Theorem  2.2 and $S$ is one of the slivers coming from the $|y| < b|x|$ region,
we have the estimate
$$\bigg|\int_S \alpha(x,y)\, \chi_E(x,y) \,\prod_{i = 1}^n |f_i(x,y)|^{\gamma_i} e^{itx + iuy}\,dx\,dy\bigg| $$
$$\leq 
C\int_S (1 + |(t,u) \cdot v| |(x,y)| + |u||(x,y) \cdot v^{\perp}| )^{-{1 \over 2}} \chi_E(x,y) \,\prod_{i = 1}^n |f_i(x,y)|^{\gamma_i}\,dx\,dy \eqno (3.1)$$
Here $v$ denotes a unit vector tangent to the sliver $S$ at the origin, and $v^{\perp}$ a normal vector; in the case where the two
boundary curves of $S$ at the origin have different tangents (i.e. $S$ is  a "wedge") then $v$ denotes the tangent to the
 boundary curve of $S$ nearest to the $x$-axis. The constant 
$C$ here depends on the function $\prod_{i = 1}^n |f_i(x,y)|^{\gamma_i}$, $E$, the application of Theorem 2.2 being used and the constant $A$ of $(1.2)-(1.3)$. If $S$ is a sliver from the  $|y| >  b|x|$ region the corresponding estimate holds with the $|u|$ factor
 replaced by $|t|$ and one replaces the $x$-axis with the $y$-axis in the above.

\noindent {\bf Proof.} We examine the integral $\int_S \alpha(x,y)\, \chi_E(x,y) \,\prod_{i = 1}^n |f_i(x,y)|^{\gamma_i} e^{itx + iuy}\,dx\,dy$ in the new coordinates after applying Theorem 2.2. The coordinate change transferring old coordinates to new 
 is either of the form $(x,y) \rightarrow (\pm x, \pm y + k(x))$, or consists of a reflection $(x,y) \rightarrow (y,x)$ followed by a 
mapping of such form. Here $k(x^N)$ is real analytic for some positive integer $N$. We will consider only the case where it is of 
the form $(x,y) \rightarrow (x,  y + k(x))$ as Lemma 3.2
for the other situations follow from this case as applied to reflected versions of $f_1(x,y),...,f_n(x,y)$ and the real analytic functions
defining the boundary curves of $E$. 

\noindent In the new coordinates, $\int_S \alpha(x,y)\, \chi_E(x,y) \,\prod_{i = 1}^n |f_i(x,y)|^{\gamma_i} e^{itx + iuy}\,dx\,dy$ becomes
$$ \int_D \alpha(x,y + k(x))\, \chi_E(x,y + k(x)) \,\prod_{i = 1}^n |f_i(x,y + k(x))|^{\gamma_i} e^{ it x +  iu y + iu k(x)}\,dx\,dy $$
Here $D$ denotes the sliver in the new coordinates (what is called $D_{il}'$ in  the notation of Theorem 2.2). Because the real analytic
functions defining the boundary curves have had their singularities resolved, those functions are comparable to monomials in the new
coordinates. In particular, they cannot have zeroes in $D$. Hence $\chi_E(x,y + k(x))$ is either identically zero or identically $1$
on $D$. Clearly we need only consider the case where it is identically $1$. In addition, since the order of  the zero of $k(x)$ at the origin
is at least one,  $\alpha(x,y + k(x))$ satisfies the estimates $(1.2)-(1.3)$ since $\alpha(x,y)$ does. So we denote $\alpha(x,y + k(x))$
by $\beta(x,y)$ and we are considering the following expression, where $\beta(x,y)$ satisfies $(1.2)-(1.3)$.
$$\int_D \beta(x,y)\,\prod_{i = 1}^n |f_i(x,y + k(x))|^{\gamma_i} e^{ it x +  iu y + iu k(x)}\,dx\,dy \eqno (3.2)$$
Note that each $f_i(x,y + k(x))$ here is comparable to a monomial in the sense of Theorem 2.2. Next, since order of  the zero 
of $k(x)$ at the origin is at least $1$, we may write $k(x) = cx + l(x)$, where $l(x)$ has a zero of order greater than one at the 
origin. Here $c$ and/or $l(x)$ may be zero. Accordingly, $(3.2)$ can be rewritten as 
$$\int_D \beta(x,y)\,\prod_{i = 1}^n |f_i(x,y + k(x))|^{\gamma_i} e^{ i(t + cu)x +  iu y + iu l(x)}\,dx\,dy \eqno (3.3)$$
We denote the expression $(3.3)$ by $I$, and we divide the integral $I$ dyadically in the $x$ and $y$ variables. Namely, for
a nonnegative smooth compactly-supported function $s(x)$ on $\R$ that vanishes on a neighborhood of $0$, we write $I = 
\sum_{jk} I_{jk}$, where
$$I_{jk} = \int_D s(2^j x)s(2^ky)\beta(x,y)\,\prod_{i = 1}^n |f_i(x,y + k(x))|^{\gamma_i} e^{ i(t + cu)x +  iu y + iu l(x)}\,dx\,dy \eqno (3.4)$$
We will apply the Van der Corput lemma (Lemma 3.1) in $(3.4)$  in the $x$ and/or $y$ direction. Adding the result over all $j$ and $k$ will give the needed bounds for $I$. We start with the $y$-direction, which it will turn out will only be needed when $l(x)$ is
identically zero. We apply the Van der Corput Lemma for first derivatives in the $y$-direction. By applying 
Corollary 2.3 for first $y$ derivatives on each monomial-like $f_i(x,y + k(x))$, we see that taking a $y$ derivative of $\prod_{i = 1}^n |f_i(x,y + k(x))|^{\gamma_i}$ introduces a factor of magnitude at most $C{1 \over y}$. By $(1.3)$ we have $|\partial_y \beta(x,y)| \leq C{ 1\over y}$, and the support condition on $s(y)$ ensures that the $y$ derivative
of the $s(2^j y)$ factor introduces a factor satisfying the same upper bounds. Thus if $Q_{jk}$ denotes the rectangle $[2^{-j-1}, 2^{-j}] \times [2^{-k-1}, 2^{-k}]$, the Van der Corput lemma
 for first derivatives leads to a bound of 
$$|I_{jk}| \leq  C \int_{Q_{jk}} {1 \over |uy|} \prod_{i = 1}^n |f_i(x,y + k(x))|^{\gamma_i}\,dx\,dy \eqno (3.5)$$
Just taking absolute values and integrating in $(3.4)$ leads to the bound
$$|I_{jk}| \leq  C \int_{Q_{jk}} \prod_{i = 1}^n |f_i(x,y + k(x))|^{\gamma_i}\,dx\,dy \eqno (3.6)$$
Thus combining $(3.5)$ and $(3.6)$ we obtain 
$$|I_{jk}| \leq  C \int_{Q_{jk}} \min\bigg(1, {1 \over |uy|}\bigg) \prod_{i = 1}^n |f_i(x,y + k(x))|^{\gamma_i}\,dx\,dy \eqno (3.7)$$
For our purposes however, we only need the weaker statement
$$|I_{jk}| \leq  C \int_{Q_{jk}} \min\bigg(1, {1 \over |uy|^{1 \over 2}}\bigg) \prod_{i = 1}^n |f_i(x,y + k(x))|^{\gamma_i}\,dx\,dy \eqno (3.8)$$
Next, in the event that $l(x)$ is not identically zero, we apply the Van der Corput lemma for second derivatives in the $x$ direction.
 Note $l''(x) \sim {l(x) \over x^2}$ on a small enough neighborhood of the origin (which we may assume we are in). By Corollary
2.3, applying an $x$ derivative to $\prod_{i = 1}^n |f_i(x,y + k(x))|^{\gamma_i}$ yields a factor of at most $C {1 \over x}$. 
This time, by $(1.3)$ we have $|\partial_x \beta(x,y)| \leq C{ 1\over x}$, and the support condition on $s(x)$ ensures that 
taking the $x$ derivative
of the $s(2^j x)$ factor incurs a factor satisfying the same upper bounds. Thus applying the Van der Corput lemma we get
$$|I_{jk}| \leq  C \int_{Q_{jk}} {1 \over |u l(x)|^{1 \over 2}} \prod_{i = 1}^n |f_i(x,y + k(x))|^{\gamma_i}\,dx\,dy \eqno (3.9)$$
(The ${1 \over x^2}$ factor one gets from taking the second derivative of $l(x)$ is exactly enough to compensate for the 
${1 \over x}$ that one 
normally gets in such applications of the Van der Corput lemma.) As in the steps from $(3.5)-(3.8)$, this leads to
$$|I_{jk}| \leq  C \int_{Q_{jk}} \min\bigg(1, {1 \over |u l(x)|^{1 \over 2}}\bigg) \prod_{i = 1}^n |f_i(x,y + k(x))|^{\gamma_i}\,dx\,dy \eqno (3.9')$$
Lastly, suppose that on the domain of integration in $(3.4)$ one has $\inf |(t + cu)x| > B \sup|ul(x)|$ (such as when $l(x)$ is
 identically zero), where the constant $B$ is large enough to ensure that if we are on a sufficiently small neighborhood of the origin, 
which we may assume, the absolute value of the
first $x$-derivative of the phase in $(3.4)$ is bounded below by ${1 \over 2} |t + cu|$. In this situation, we 
may apply the Van der Corput lemma for first derivatives in the $x$ direction. This time we obtain a bound of 
$$|I_{jk}| \leq  C \int_{Q_{jk}} {1 \over |(t + cu)x|} \prod_{i = 1}^n |f_i(x,y + k(x))|^{\gamma_i}\,dx\,dy \eqno (3.10)$$ 
Like in the steps from $(3.5)-(3.8)$ this implies that
$$|I_{jk}| \leq  C \int_{Q_{jk}} \min\bigg(1, {1 \over |(t + cu)x|^{1 \over 2}}\bigg) \prod_{i = 1}^n |f_i(x,y + k(x))|^{\gamma_i}\,dx\,dy \eqno (3.10')$$
Combining $(3.9')$ and $(3.10')$, we have for all $(j,k)$ that
$$|I_{jk}| \leq  C \int_{Q_{jk}} \min\bigg(1,  {1 \over |u l(x)|^{1 \over 2}}, {1 \over |(t + cu)x|^{1 \over 2}}\bigg) \prod_{i = 1}^n |f_i(x,y + k(x))|^{\gamma_i}\,dx\,dy \eqno (3.11)$$
Finally, combining with $(3.8)$, we see that for each $(j,k)$ we have
$$|I_{jk}| \leq  C \int_{Q_{jk}} \min\bigg(1,  {1 \over |uy|^{1 \over 2}}, {1 \over |u l(x)|^{1 \over 2}}, {1 \over |(t + cu)x|^{1 \over 2}}\bigg) \prod_{i = 1}^n |f_i(x,y + k(x))|^{\gamma_i}\,dx\,dy \eqno (3.12)$$
In view of the shape of $D$ as given by Theorem 2.2 (where it is called $D_{il}'$), summing $(3.12)$ over all $(j,k)$ leads to
$$|I| \leq C \int_D \min\bigg(1,  {1 \over |uy|^{1 \over 2}}, {1 \over |u l(x)|^{1 \over 2}}, {1 \over |(t + cu)x|^{1 \over 2}}\bigg) \prod_{i = 1}^n |f_i(x,y + k(x))|^{\gamma_i}\,dx\,dy \eqno (3.13)$$
We are now in a position to prove $(3.1)$. First suppose $k(x)$ is identically zero. Then $l(x)$ is identically zero and $c=0$, and 
$(3.13)$ becomes
$$|I| \leq C \int_D \min\bigg(1,  {1 \over |uy|^{1 \over 2}}, {1 \over |tx|^{1 \over 2}}\bigg) \prod_{i = 1}^n |f_i(x,y)|^{\gamma_i}\,dx\,dy \eqno (3.14)$$
By the form of $D$ given by Theorem 2.2, one has $v = (1,0)$  in $(3.1)$ when $k(x)$ is identically zero (see the discussion at the end of the proof for the case when $D$ is a wedge.) Therefore equation $(3.14)$ is equivalent to $(3.1)$ and we are done. So we move to the case where $k(x)$ is not identically zero. Then $(3.13)$ implies
$$|I| \leq C \int_D \min\bigg(1,  {1 \over |u l(x)|^{1 \over 2}}, {1 \over |(t + cu)x|^{1 \over 2}}\bigg) \prod_{i = 1}^n |f_i(x,y + k(x))|^{\gamma_i}\,dx\,dy \eqno (3.15)$$
Doing the variable change $(x,y) \rightarrow (x, y - k(x))$ to turn the sliver back into its original coordinates, $(3.15)$ becomes
$$|I| \leq C \int_S \min\bigg(1,  {1 \over |u l(x)|^{1 \over 2}}, {1 \over |(t + cu)x|^{1 \over 2}}\bigg) \prod_{i = 1}^n |f_i(x,y)|^{\gamma_i}\,dx\,dy \eqno (3.16)$$
The quantity $|(t + cu)x| = |(t,u)\cdot (1,c)||x|$ is within a bounded factor of $|(t,u) \cdot (1,c)| |(x,y)|$ since we are assuming the sliver $S$ is from the $|y| < b|x|$ region in Theorem 2.2. Also, note that $(1,c)$ is tangent to the sliver. Hence
 $|(t,u)\cdot (1,c)||(x,y)|$ is within a bounded factor of $|(t,u) \cdot v||(x,y)|$, where  $v$ is a unit tangent vector as in the statement
of Lemma 3.2. On the other hand,
the quantity $l(x)$ is the vertical drop between $(x,y)$ and the line with direction $v$ through the origin, and since the sliver is in
the $|y| < b|x|$ region this vertical drop is within a bounded factor of the distance from $(x,y)$ to this line, which is given by
$|(x,y) \cdot v^{\perp}|$. Hence $|u l(x)|$ is within a bounded factor of $|u| |(x,y) \cdot v^{\perp}|$. Thus $(3.16)$ implies
$$|I| \leq C \int_S \min\bigg(1,  {1 \over (|u| |(x,y) \cdot v^{\perp}|)^{1 \over 2}}, {1 \over (|(t,u) \cdot v||(x,y)|)^{1 \over 2}}\bigg) \prod_{i = 1}^n |f_i(x,y)|^{\gamma_i}\,dx\,dy \eqno (3.17)$$
This is equivalent to $(3.1)$ as desired. 

As for the statement in Lemma 3.2 concerning which $v$ to choose  when $S$ is a wedge-shaped region with two tangent lines at 
the origin, such an $S$ can arise in two ways. We focus on the wedges where $x > 0$ and $|y| < bx$ as the other cases 
are very similar.  One way for such a wedge to arise occurs at the beginning of the resolution process of Theorem 2.1 when $S$ is of the form $\{(x,y): 0 < x < a, \,hx^m < y < Hx^M\}$, for $h, H \geq 0$, $m > M$ or $\{(x,y): 0 < x < a,\, hx^m >  y >  Hx^M\}$  for
$h, H \leq 0$, $m > M$.
In these cases $ k(x)$ is always identically zero, so the correct tangent line to choose for $S$ is the one closest to the $x$-axis.
The other way such an $S$ can arise is again early in the resolution process when $S$ is of the form 
$\{(x,y): 0 < x < a, hx  < y < Hx\}$ for some $h \neq H$ 
and the resolution process is such that $k(x)$ takes the $x$-axis to the nearer boundary curve of $S$ via a map of the form
$(x,y) \rightarrow (x, \pm y + cx)$ for an appropriate $c$. Once again the correct 
boundary curve of $S$ to choose is the one nearest the $x$-axis. This completes the proof of Lemma 3.2.

\noindent {\bf Lemma 3.3.} Suppose we are not in the trivial situation where $E$ contains a neighborhood of the origin and each $f_i(0,0) \neq 0$. Then a sufficiently small disk $B$ centered at the origin can be 
written in the form $B = \cup_{i=1}^M B_i$, where each $B_i$ is a wedge bounded by lines through the origin and the boundary of 
$B$, such that each $B_i$ is of one of the following two forms.

\noindent {\bf 1)} Let $n_i$ denote the order of the zero of $f_i$ at the origin. Then on the first type of wedge, for
some positive constants $c_i$ and $c_i'$, $f_i(x,y)$ satisfies
$$c_i(x^2 + y^2)^{n_i \over 2} < |f_i(x,y)| < C_i(x^2 + y^2)^{n_i \over 2} \eqno (3.18)$$
Furthermore,  the boundary curves of $E$ do not intersect the closure of $B_i$ and one has
$$\bigg|\int_{B_i}\alpha(x,y)\, \chi_E(x,y) \,\prod_{i = 1}^n |f_i(x,y)|^{\gamma_i} e^{itx + iuy}\,dx\,dy\bigg| $$
$$\leq C\int_{B_i} (1 + |tx| + |uy| )^{-1} \chi_E(x,y) \,\prod_{i = 1}^n |f_i(x,y)|^{\gamma_i}\,dx\,dy \eqno (3.19)$$
\noindent {\bf 2)} Let $F_i(x,y)$ denote the sum of the terms of $f_i(x,y)$'s Taylor expansion of lowest degree.
If $B_i$ is the second type of wedge, there is a line $l_i$ through the origin intersecting $B_i$ that is either 
part of the zero set of one of the $F_i(x,y)$ or tangent to one of the boundary curves of $E$ at the origin. Furthermore, if
 $v$ denotes a unit vector in the direction of $l_i$ then we have
$$\bigg|\int_{B_i} \alpha(x,y)\, \chi_E(x,y) \,\prod_{i = 1}^n |f_i(x,y)|^{\gamma_i} e^{itx + iuy}\,dx\,dy\bigg| $$
$$\leq 
C\int_{B_i} (1 + |(t,u) \cdot v| |(x,y)| + |(t,u) \cdot v^{\perp}||(x,y) \cdot v^{\perp}| )^{-{1 \over 2}} \chi_E(x,y) \,\prod_{i = 1}^n |f_i(x,y)|^{\gamma_i}\,dx\,dy \eqno (3.20)$$
\noindent {\bf Proof.}  We apply Theorem 2.2 to all of the $f_i(x,y)$ as well as  real analytic functions whose zero
set contains the boundary of $E$. We let the first type of $B_i$  be certain wedges which can be described in terms of the resolution of singularities process of Theorem 2.2 as follows.

 Let $F(x,y)$ denote the sum of terms of lowest degree
of the Taylor expansion at the origin of the product of functions whose zero set is being resolved. 
At the beginning of the resolution
 of singularities process of Theorem 2.1, one has a collection of wedges associated to the edge of slope $-1$ of the Newton polygon of the 
product of functions being resolved, that are away from the zeroes of the $F(x,y)$. These wedges are bounded by two lines
through the origin and a vertical or horizontal line. 

Since the zeroes of 
$F(x,y)$ are the union of
the zeroes of the $G_i(x,y)$, where $G_i(x,y)$ denotes the sum of terms of one of the functions in the product,
these wedges are away from the zeroes of any $G_i(x,y)$ as well. We declare that any intersection of one of these wedges
 with the disk $B$ is a $B_i$ of the first type of in Lemma 3.2.
Because they are away from the zeroes of any $G_i(x,y)$, equation  $(3.18)$ holds. Furthermore one could have taken
$k(x)$ to be zero for these wedges, since no resolution of singularities is needed. Equation $(3.19)$ is therefore a consequence of $(3.7)$ and $(3.10)$, summed over all $j$ and $k$.

The complement of the union of the $B_i$ above is, modulo boundaries, a finite union of disjoint wedges. By the constructions of
Theorem 2.2, each
sliver $S$ that is not in one of the wedges $B_i$ above is contained in one of these new wedges. By
construction, each such wedge contains exactly one line through the origin which is in the zero set of $F(x,y)$. Since this zero
set is the union of the zero sets of the $G_i(x,y)$, the line in question is either a zero set of an $F_i(x,y)$ coming from an 
$f_i(x,y)$, a tangent line at the origin to a boundary curve of $E$, or a tangent line at the origin  to one of the other curves which are in the
zero set of the real analytic functions whose zero sets contains the boundary curves to $E$, but which is not also one of the
earlier tangent lines. If the line is of the last variety,
we let this wedge be a $B_i$ of the first kind, and $(3.18)-(3.19)$ holds exactly as before. All other wedges are declared to be 
wedges of the second kind.

Thus in order to prove part 2 of this lemma
it suffices to show $(3.20)$ for the second kind of wedge, 
 where $l_i$ is the line through the origin contained in the closure of $B_i$ which is in the zero set of  $F(x,y)$.

The vectors denoted by $v$ in Lemma 3.2 are of the form $(1,c)$, where the coordinate shift $(x,y) \rightarrow (x, y + k(x))$ 
satisfies $k(x) = cx + $ higher order terms. The resolution of singularities process of Theorem 2.2 is such that the line
 $y = cx$ is contained in the zero set of $F(x,y)$, and all slivers $S$ are contained in some $B_i$. Thus the $v$ of Lemma 3.2 is of the type needed in part b) of this lemma.
We now add $(3.1)$ over all slivers $S$ contained in a given $B_i$ of the second type and obtain
$$\bigg|\int_{B_i} \alpha(x,y)\, \chi_E(x,y) \,\prod_{i = 1}^n |f_i(x,y)|^{\gamma_i} e^{itx + iuy}\,dx\,dy\bigg| $$
$$\leq 
C\int_{B_i} (1 + |(t,u) \cdot v| |(x,y)| + |u||(x,y) \cdot v^{\perp}|)^{-{1 \over 2}}\chi_E(x,y) \,\prod_{i = 1}^n |f_i(x,y)|^{\gamma_i}\,dx\,dy \eqno (3.21)$$
This is almost the same as $(3.20)$. The one difference is that instead of having a $|(t,u) \cdot v^{\perp}|$ factor as in $(3.20)$ we
have a $ |u|$ factor. Suppose we could show that for some constant $e$ the following inequality holds on $B_i$.
 $$|(t,u) \cdot v| |(x,y)| \geq e|(t,u) \cdot v^{\perp}||(x,y) \cdot v^{\perp}| \eqno (3.22)$$ 
Then the  $|(t,u) \cdot v| |(x,y)| $ term alone is enough for $(3.21)$ to imply $(3.20)$. 
 This would only not hold if
$(t,u)$ is nearly in the $v^{\perp}$ direction. In this case $|(t,u) \cdot v^{\perp}|$  is of comparable magnitude to $|(t,u)|$.
Because $v$ is in the direction of  $(1,c)$ for fixed $c$, there's a $M$ such that if $|t| > M|u|$ then $(3.22)$ holds.
Otherwise, $|(t,u)|$ is of comparable magnitude to $|u|$, so $|(t,u) \cdot v^{\perp}|$  is also of comparable magnitude to
$|u|$. In this case
$(3.21)$ once again implies $(3.20)$ as needed. This concludes the proof of Lemma 3.3.

\noindent Note that since on the wedge $B_i$ we have $|(x,y)| \sim |(x,y) \cdot v|$, one can write $(3.20)$ in the symmetric form
$$\bigg|\int_{B_i} \alpha(x,y)\, \chi_E(x,y) \,\prod_{i = 1}^n |f_i(x,y)|^{\gamma_i} e^{itx + iuy}\,dx\,dy\bigg| $$
$$\leq 
C\int_{B_i} (1 + |(t,u) \cdot v| |(x,y) \cdot v| + |(t,u) \cdot v^{\perp}||(x,y) \cdot v^{\perp}|)^{-{1 \over 2}} \chi_E(x,y) \,\prod_{i = 1}^n |f_i(x,y)|^{\gamma_i}\,dx\,dy \eqno (3.23)$$
In this form it is readily apparent how the estimate $(3.20)$ is independent of the resolution of singularities process being used. 

Next,
we give the following corollary to Lemma 3.3 which we will need for the  paper [G3].

\noindent {\bf Corollary 3.4.}  Suppose
 that $E$ is a disk centered at the origin, but we are not in the trivial situation where each $f_i(0,0) \neq 0$. Let
 $p_1 \neq  p_2$, and let $V$ be any of the four wedges with vertex $(0,0)$ formed by the lines $y = p_1x$ and $y = p_2x$.
Suppose that each $F_i(x,y)$ has no zeroes on set $cl(V) \cap (\R - \{0\})^2$.  
Then if the function $\alpha(x,y)$ in $(1.1)$ is supported on the disk $B$ where Lemma 3.3 applies, then we have 
the following simplified version of $(3.1)$.
$$\bigg|\int_{B \cap V} \alpha(x,y)\,\prod_{i = 1}^n |f_i(x,y)|^{\gamma_i} e^{itx + iuy}dx\,dy\bigg|  \leq C\int_{B \cap V} (1 + |tx| + |uy|)^{-{1 \over 2}} \,\prod_{i = 1}^n |f_i(x,y)|^{\gamma_i}dx\,dy \eqno (3.24a)$$
If each $F_i(x,y)$ has no zeroes on all of $(\R - \{0\})^2$, one has 
$$\bigg|\int_{B} \alpha(x,y)\,\prod_{i = 1}^n |f_i(x,y)|^{\gamma_i} e^{itx + iuy}\,dx\,dy\bigg| \leq  C\int_{B} (1 + |tx| + |uy|)^{-{1 \over 2}} \,\prod_{i = 1}^n |f_i(x,y)|^{\gamma_i}\,dx\,dy \eqno (3.24b)$$
\noindent {\bf Proof.}  The second part follows immediately from the first, as in the setting of the second part one can write 
$\R^2$ as the union of four $V$ on which $(3.24a)$ applies. Then $(3.24b)$ follows by addition. As for part a), one applies
Lemma 3.3 to  $f_1(x,y),...,f_n(x,y)$. As long as the $B_i$ of the second type were chosen to be narrow enough, each $B_i \cap V$ for a $B_i$ of the second type will be empty unless the line $l_i$ is the $x$ or $y$ axis; the other $l_i$ are zeroes of some $F_i(x,y)$ which lie outside of $cl(V)$. In this situation we can define the $B_i$ of the first type so that $B \cap V$ is a union of
some $B_i \cap V$ where each $B_i$ of either of the first type or of the second type with $v = (1,0)$ or $(0,1)$. Then adding $(3.19)$ or $(3.20)$ over all $B_i$ gives the corollary.

\noindent The the next lemma will help go from Lemma 3.3 to the estimates of Theorems 1.3-1.5.

\noindent {\bf Lemma 3.5.} Let $B_i$ be one of the domains of part 2 of  Lemma 3.3, and $v$ a unit vector in the direction of the
associated line $l_i$.

\noindent {\bf a)} Let $(\epsilon,d)$ be as in $(1.5)$. For any $a > 0$, let $F_a=\{(t,u): |(t,u) \cdot v| > 
a |(t,u) \cdot v^{\perp}|\}$. Then if $\epsilon <  {1 \over 2}$ there is a constant $C_a$ such that for $(t,u) \in F_a$ one has an estimate
$$\bigg|\int_{B_i} \alpha(x,y)\, \chi_E(x,y) \,\prod_{i = 1}^n |f_i(x,y)|^{\gamma_i} e^{itx + iuy}\,dx\,dy\bigg|
 \leq C_a(2 + |(t,u)|)^{-\epsilon}( \ln (2 + |(t,u)|))^d \eqno (3.25)$$
If $\epsilon = {1 \over 2}$ one gets the estimate obtained by replacing $d$ by $d + 1$ in $(3.25)$, and if $\epsilon > {1 \over 2}$
one has $(2 + |(t,u)|)^{-{1 \over 2}}$ in place of $(2 + |(t,u)|)^{-\epsilon}( \ln (2 + |(t,u)|))^d $.

\noindent {\bf b)} \noindent {\bf b)} For $w = v^{\perp}$, let $(\delta_w, e_w)$ be as in $(1.6)$, and let $G =  \{(t,u): |(t,u) \cdot v| < |(t,u) \cdot v^{\perp}|\}$.
Then if $\delta_w < {1 \over 2}$ there is a constant $D$ such that for $(t,u) \in G$ one has an estimate
$$\bigg|\int_{B_i} \alpha(x,y)\, \chi_E(x,y) \,\prod_{i = 1}^n |f_i(x,y)|^{\gamma_i} e^{itx + iuy}\,dx\,dy\bigg|
 \leq D(2 + |(t,u)|)^{-\delta_w}( \ln (2 + |(t,u)|))^{e_w} \eqno (3.26)$$
If $\delta_w = {1 \over 2}$ one gets the estimate obtained by replacing $e_w$ by $e_w + 1$, and if $\delta_w > {1 \over 2}$
one has $(2 + |(t,u)|)^{-{1 \over 2}}$ in place of $(2 + |(t,u)|)^{-\delta_w}( \ln (2 + |(t,u)|))^{e_w}$.
 
\noindent {\bf Proof.} We start with part a). We can assume that $|(t,u)| > 4$ say, since the case where $|(t,u)| \leq 4$ is immediate.
On the domain $F_a$, there is a constant $a'$ such that $|(t,u) \cdot v| > a' |(t,u)|$. Thus
ignoring the  $|(t,u) \cdot v^{\perp}||(x,y) \cdot v^{\perp}|$ term in $(3.20)$, we see that on $F_a$, $(3.20)$ implies that
$$\bigg|\int_{B_i} \alpha(x,y)\, \chi_E(x,y) \,\prod_{i = 1}^n |f_i(x,y)|^{\gamma_i} e^{itx + iuy}\,dx\,dy\bigg| $$
$$\leq 
C \int_{B_i} (1 + |(t,u)| |(x,y)| )^{-{1 \over 2}} \chi_E(x,y) \,\prod_{i = 1}^n |f_i(x,y)|^{\gamma_i}\,dx\,dy \eqno (3.27)$$
We divide the integral $(3.27)$ into $|(x,y)| < |(t,u)|^{-1}$ and $|(x,y)| > |(t,u)|^{-1}$ parts. The integral over the first part is
$$C \int_{\{(x,y) \in B_i:\,|(x,y)| < |(t,u)|^{-1}\}} \chi_E(x,y) \,\prod_{i = 1}^n |f_i(x,y)|^{\gamma_i}\,dx\,dy \eqno (3.28)$$
By Lemma 1.1 this is bounded by the desired bound of $C_a(2 + |(t,u)|)^{-\epsilon}( \ln (2 + |(t,u)|))^d$. For the $|(x,y)| > |(t,u)|^{-1}$ part, we divide the integral dyadically in $|(x,y)|$, obtaining a bound of  
$$C  \sum_{j = 1}^{\lfloor \log_2 |(t,u)| - 1 \rfloor} 2^{-{j \over 2}} \int_{\{(x,y) \in B_i:\, 2^j|(t,u)|^{-1} \leq |(x,y)| < 2^{j+1}|(t,u)|^{-1}\}} \chi_E(x,y) \,\prod_{i = 1}^n |f_i(x,y)|^{\gamma_i}\,dx\,dy \eqno (3.29) $$
We can use $\lfloor \log_2 |(t,u)| - 1 \rfloor$ here because our domain is contained in a small neighborhood of the origin. Inserting $(1.5)$
into the above provides a bound of
$$C  \sum_{j = 1}^{ \lfloor \log_2 |(t,u)| - 1 \rfloor} 2^{-{j \over 2}} (2^{j}|(t,u)|^{-1})^{\epsilon}|\ln (2^{j}|(t,u)|^{-1})|^d \eqno (3.30)$$
If $ \epsilon < {1 \over 2}$ the summands decrease exponentially and we obtain the bound given by the first term, namely
$ C |(t,u)|^{-\epsilon}| \ln |(t,u)||^d$. If $\epsilon > {1 \over 2}$ the summands increase  exponentially and we obtain the bound 
given by the last term, namely $|(t,u)|^{-{1 \over 2}}$. When $\epsilon = {1 \over 2}$, we get a summation of powers of 
$|\ln (2^{j}|(t,u)|^{-1})|^d$ which results in  a bound of $ C |(t,u)|^{-{1 \over 2}}| \ln |(t,u)||^{d+1}$. This is
equivalent to the statement of Lemma 3.5a) since we can assume that $|(t,u)| > 4$. This completes the proof of part a).

The proof of part b) is rather similar. As before, it suffices to assume $|(t,u)| > 4$. In the domain $G$,  since $|(t,u) \cdot v| < |(t,u) \cdot v^{\perp}|$, one has that
$|(t,u) \cdot v^{\perp}| \sim |(t,u)|$ and the right-hand side of $(3.20)$ can be bounded by
$$C\int_{B_i} (1 + |(t,u)||(x,y) \cdot v^{\perp}|)^{-{1 \over 2}}\chi_E(x,y) \,\prod_{i = 1}^n |f_i(x,y)|^{\gamma_i}\,dx\,dy \eqno (3.31)$$
We break the integral in $(3.31)$ into $|(x,y) \cdot v^{\perp}| < |(t,u)|^{-1}$ and $|(x,y) \cdot v^{\perp}| \geq |(t,u)|^{-1}$
parts. By Lemma 1.2, the first part is bounded by $C |(t,u)|^{-\delta_w}(\ln |(t,u)|)^{e_w}$. We again break up the second 
integral dyadically, obtaining a bound of 
$$\sum_{j = 1}^{\lfloor \log_2 |(t,u)| - 1 \rfloor} 2^{-{j \over 2}}\int_{\{(x,y) \in B_i:\,2^j|(t,u)|^{-1} \leq |(x,y) \cdot v^{\perp}| \leq 2^{j+1}|(t,u)|^{-1}\}}\chi_E(x,y) \,\prod_{i = 1}^n |f_i(x,y)|^{\gamma_i}\,dx\,dy \eqno (3.32) $$
Using Lemma 1.2, this leads to a bound of
$$C  \sum_{j = 1}^{ \lfloor \log_2 |(t,u)| - 1 \rfloor} 2^{-{j \over 2}} (2^{j}|(t,u)|^{-1})^{\delta_w}|\ln (2^{j}|(t,u)|^{-1})|^{e_w} \eqno (3.33)$$
One then argues as after $(3.30)$, and we see that we get a bound of $ C |(t,u)|^{-\delta_w}| \ln |(t,u)||^{e_w}$ when 
$\delta_w < {1 \over 2}$, a bound of $C|(t,u)|^{-{1 \over 2}}$ when $\delta_w > {1 \over 2}$, and 
$ C |(t,u)|^{-\delta_w}| \ln |(t,u)||^{e_w+1}$  when $\delta_w = {1 \over 2}$. This completes the proof of Lemma 3.5.

\noindent For the $B_i$ of part a) of Lemma 3.3, we have estimates at least as strong as those that hold for the $B_i$ of part b) of
Lemma 3.3:

\noindent {\bf Lemma 3.6.} If $B_i$ is a wedge from part a) of Lemma 3.3 and $\epsilon < 1$, then for all $(t,u)$ one has an estimate
$$\bigg|\int_{B_i} \alpha(x,y)\, \chi_E(x,y) \,\prod_{i = 1}^n |f_i(x,y)|^{\gamma_i} e^{itx + iuy}\,dx\,dy\bigg|
 \leq C_a(2 + |(t,u)|)^{-\epsilon}( \ln (2 + |(t,u)|))^d \eqno (3.34)$$
If $\epsilon = 1$ one gets the estimate obtained by replacing $d$ by $d + 1$ in $(3.34)$, and if $\epsilon  > 1$
one has $(2 + |(t,u)|)^{-1}$ in place of $(2 + |(t,u)|)^{-\epsilon}( \ln (2 + |(t,u)|))^d $.

\noindent {\bf Proof.} Equation $(3.19)$ implies 
$$\bigg|\int_{B_i}\alpha(x,y)\, \chi_E(x,y) \,\prod_{i = 1}^n |f_i(x,y)|^{\gamma_i} e^{itx + iuy}\,dx\,dy\bigg| $$
$$\leq C\int_{B_i} (1 + |(t,u)||(x,y)|)^{-1} \chi_E(x,y) \,\prod_{i = 1}^n |f_i(x,y)|^{\gamma_i}\,dx\,dy \eqno (3.35)$$
This is exactly $(3.27)$ with the exponent $-{1 \over 2}$ replaced by $-1$. The steps from $(3.27)$ through the paragraph after
$(3.30)$ lead to the statement of this lemma, with this modification due to the new exponent.

\noindent {\bf Proofs of Theorem 1.3-1.5.}

We now are in a position to prove Theorems 1.3-1.5. If $E$ contains a neighborhood of $(0,0)$ and each $f_i(0,0) \neq 0$ the
results are easy, so we assume that this is not the case. We start with Theorems 1.3a) and 1.4a). First of all, note that the directions
of the lines $l_i$ in Theorems 1.3 and 1.4 are exactly the directions $v$ in Lemma 3.3 corresponding to the domains $B_i$ of the 
second type. If $l$ is not perpendicular to one of these directions,
then for each $B_i$ of the second type, the set $l_H$ of Theorems 1.3 and 1.4 is contained in one of the sets $F_a$ of Lemma 3.5, except for 
an inconsequential part near the origin. Thus for $(t,u) \in l_H$, the estimates of part a) of Lemma 3.5 hold, 
possibly with different constants.  By Lemma
3.6 they also hold for $(t,u) \in l_H$ when $B_i$ is of the type of the first part of Lemma 3.3. Hence they hold for all $B_i$. Adding 
this over all $i$  therefore results in the bounds of Lemma 3.5a) whenever $(t,u) \in l_H$, giving
an estimate $|K_{\phi,x_0,y_0}(t,u)| \leq C(2 + |(t,u)|)^{-\epsilon}(\ln (2 + |(t,u)|))^d$ on $l_H$ when
$l$ is not in one of the $v^{\perp}$ directions. In section 4, we will show the best possible power of $(2 + |(t,u)|)$ that 
can appear in such an estimate is $-\delta_v$. Since $\delta_v \leq \epsilon$ for any direction, we have 
$\epsilon = \delta_v$ here as needed in Theorems 1.3a) and 1.4a). This provides the estimates of Theorems 1.3a) and 1.4a) 
when $l$ is not perpendicular to the direction of a $v$ corresponding to the second type of $B_i$.

 If $l$ is in the $v^{\perp}$ direction for some $B_{i_0}$ of the second type, when $(t,u) \in l_H$ one has the same estimates 
given by Lemma 3.5a) for $i \neq i_0$, for the same reasons as before. For $B_{i_0}$ however, the arguments above will not work
since $l_H$ is a subset of no $F_a$. However, besides an inconsequential region with small $|(t,u)|$, 
$l_H$ is a subset of the set called $G$ in Lemma 3.5b). So for $(t,u) \in l_H$ one has the possibly weaker bounds of Lemma 
3.5b)  in place of those of Lemma 3.5a). Adding this to the estimates over the other $B_i$ therefore gives the overall bound given by 
part b) of Lemma 3.5. This gives the statements of Theorems 1.3a) and 1.4a) when $l$ is one of the $v^{\perp}$ directions.

Moving on to the overall decay rates of Theorems 1.3b, 1.4b, and 1.5, for a 
given $B_i$ of the second type, each $(t,u)$ is either in a $F_a$ of Lemma 3.5a) or a set  $G$ of Lemma 3.5b). (The 
exact value of $a$ is not important for our purposes.) Thus the quantity $|\int_{B_i}\alpha(x,y)\, \chi_E(x,y) \,\prod_{i = 1}^n |f_i(x,y)|^{\gamma_i} e^{itx + iuy}\,dx\,dy|$ is always bounded by the worse of the two estimates, namely the estimate 
given by Lemma 3.5b). The corresponding estimate for the $B_i$ 
of the first type, provided by Lemma 3.6, is at least as good as this, so adding over all $i$ we have that
$|\alpha(x,y)\, \chi_E(x,y) \,\prod_{i = 1}^n |f_i(x,y)|^{\gamma_i} e^{itx + iuy}\,dx\,dy|$ is bounded by the worst 
over all $i$ of the estimates given by Lemma 3.5b). These are exactly the overall decay rates of Theorems 1.3b, 1.4b, and 1.5. 
This completes the proofs of Theorems 1.3, 1.4, and 1.5.

\noindent {\bf Proof of Theorem 1.7.}

 It suffices to prove that given any $(x_0,y_0)$ in the support of $m(x,y)$ there is a
neighborhood $U$ of $(x_0,y_0)$ such that if the cutoff function $\phi(x,y)$ is supported in $U$, then each term of 
$(1.4)$ satisfies the bounds stipulated in Theorem 1.7 as the estimates for $K(t,u)$ then follow by addition. Let
$V$ denote the set of all $(t,u)$ within angle ${\pi \over 8}$ of the lines $y = x$ or $y = -x$. We will prove the estimates for
$(t,u)$ in the closure of $V$. The result for the remaining $(t,u)$ will follow by applying the resolution of singularities theorem 
in the coordinates obtained after rotating by $45$ degrees.

Let $S$ be any of the slivers arising from the application of Theorem 2.2 as in the previous lemmas. We will examine the contribution
to $F(t,u)$ coming from $S$ and see that it satisfies the needed bounds, so that adding over all slivers 
gives the desired estimates. We focus as before on slivers coming from the region $|y| < b|x|$ as the other slivers are treated in
an entirely analogous fashion. 

For a given sliver $S$, we move into the new coordinates and use $(3.7)$. We add $(3.7)$ over all $j$ and $k$ and call
the result $F_S(t,u)$. In view of the 
shape of the domains provided by Theorem 2.2, we get that
$$|F_S(t,u)| \leq  C \int_D \min\bigg(1, {1 \over |uy|}\bigg) \prod_{i = 1}^n |f_i(x,y + k(x))|^{\gamma_i}\,dx\,dy \eqno (3.36)$$
Here $D$ denotes the sliver $S$ in the new coordinates. Since the coordinate changes have Jacobian 1, the function $\prod_{i = 1}^n |f_i(x,y + k(x))|^{\gamma_i}$ is in $L^p(D)$, where $p$ is as in the statement of Theorem 1.7. We apply H\"older's inequality 
in $(3.36)$, obtaining
$$|F_S(t,u)| \leq C \bigg(\int_D \min\bigg(1, {1 \over |uy|^{p'}}\bigg)\bigg)^{1 \over p'} \eqno (3.37)$$
$$\leq C\bigg(\int_{[0,1] \times [0,1]} \min\bigg(1, {1 \over |uy|^{p'}}\bigg)\bigg)^{1 \over p'} \eqno (3.38)$$
Doing the integral in $(3.38)$ gives  $C\min(1, |u|^{-{1 \over p'}})$ when $p < \infty$, and
$C\min(1, |u|^{-1}\ln|u|)$ if $p = \infty$. Since we are considering $(t,u)$ within angle at most
${\pi \over 8}$ of the lines $y = x$ or $y = -x$, this gives the desired bound of  $C\min(1, |(t,u)|^{-{1 \over p'}})$ when 
$p < \infty$ and  $C\min(1, |(t,u)|^{-1}\ln|(t,u)|)$ if $p = \infty$. We add this
over all slivers $S$ and we are done.

\noindent {\heading 4. Proofs of sharpness statements.}

 The directions of part a) of Theorems 1.3 and 1.4 are given by the same directions, and we will prove 
sharpness of both simultaneously. So suppose $(1.7a)$ or $(1.9)$ holds with the exponent $\delta_v$ replaced by some $\delta > \delta_v$; we will arrive at a contradiction. Since the estimate is to hold on the whole strip $l_H$, it must hold on the ray with 
direction $v$ emanating from the origin. In other words, we have the following estimate in the $s$ variable.
$$|K_{\phi,x_0,y_0}(s v)| \leq C(2 + |s|)^{-\delta} \eqno (4.1)$$
Let $\psi(x)$ be a smooth function on 
$\R$ whose Fourier transform is a compactly supported nonnegative function equal to $1$ on a neighborhood of the origin. 
Let $0 < \eta  < \delta$ such that $\eta + \delta_v < \delta$.  For a
 large $L$ we look at
$$I_L = \int K_{\phi,x_0,y_0}(s v) \psi\bigg({s \over L}\bigg)|s|^{\delta - 1 - \eta}\,ds \eqno (4.2)$$
Inserting $(4.1)$ into $(4.2)$ gives that for all $L$ we have
$$I_L \leq C \int (2 + |s|)^{-\delta}  |s|^{\delta - 1 - \eta}\psi\bigg({s \over L}\bigg)\,ds \eqno (4.3)$$
Because $\eta > 0$, the integrand in $(4.3)$  is integrable for large $|s|$, and because $\eta  < \delta$ the 
integrand in $(4.3)$ is integrable for small $|s|$. Hence the $I_L$ are uniformly bounded in $L$. On the other hand, $I_L$ is given by
$$I_L = \int \phi(x_0 + x, y_0 + y)m(x_0 + x, y_0 + y) e^{is((x,y) \cdot v)}  |s|^{\delta - 1 - \eta}\psi\bigg({s \over L}\bigg)\,ds 
\,dx\,dy \eqno (4.4)$$
Performing the $s$ integral in $(4.4)$ leads to 
$$I_L = \int \phi(x_0 + x, y_0 + y)m(x_0 + x, y_0 + y)  L^{\delta - \eta} \xi( L [(x,y) \cdot v])\,dx\,dy \eqno (4.5)$$
Here $\xi(x)$ is the Fourier transform of $\psi(s)|s|^{\delta - 1 - \eta}$. Since the Fourier transform of 
$\psi(s)$ is nonnegative and compactly supported, and the Fourier transform of $|s|^{\delta - 1 - \eta}$ is of the form
 $c|x|^{\eta - \delta }$, we have that $\xi(x)$ is of the form $c \tilde{\xi}(s)$ where  $\tilde{\xi}(s)$ is nonnegative and
decays as $|s|^{\eta - \delta}$ as $|s| \rightarrow \infty$. Since we are assuming $\alpha(x,y)$ in $(1.1)$  is bounded below by a positive value on some neighborhood of the origin, as long as the support of $\phi(x_0 + x, y_0 + y)$ is contained in this neighborhood, 
the $\phi(x_0 + x, y_0 + y)m(x_0 + x, y_0 + y) $ factor in $(4.5)$ is nonnegative and we can 
we can rewrite $(4.5)$ as 
$$|I_L| = |c|\int \phi(x_0 + x, y_0 + y)m(x_0 + x, y_0 + y)  L^{\delta - \eta} \tilde{\xi}( L [(x,y) \cdot v])\,dx\,dy \eqno (4.5')$$
Letting $N$ be a neighborhood of the origin on which $\phi(x_0 + x, y_0 + y)$ and $\alpha(x,y)$ are both bounded below by a positive number, there is a constant $C$ such that
$$|I_L| \geq  C \int_N g(x,y)  L^{\delta - \eta} \tilde{\xi}( L [(x,y) \cdot v])\,dx\,dy \eqno (4.6)$$
As a result, for any $r > 0$ we have
$$\sup_L |I_L|  \geq  C \int_{\{(x,y) \in N:\, r < |(x,y) \cdot v| < 2r\}} g(x,y)  L^{\delta - \eta} \tilde{\xi}( L [(x,y) \cdot v])\,dx\,dy \eqno (4.7)$$
Note that the left-hand side of $(4.7)$ is finite. Since $\tilde{\xi}(s)$ is nonnegative and
decays as $|s|^{\eta - \delta}$ as $|s| \rightarrow \infty$, if we take the limit as $L \rightarrow \infty$ in the right-hand side of $(4.7)$
we obtain
$$\sup_L |I_L|  \geq  C \int_{\{(x,y) \in N:\, r < |(x,y) \cdot v| < 2r\}} g(x,y) |(x,y) \cdot v|^{\eta - \delta}\,dx\,dy \eqno (4.8)$$
As result  we have
$$\sup_r r^{\eta - \delta} \int_{\{(x,y) \in N:\, r < |(x,y) \cdot v| < 2r\}} g(x,y) \,dx\,dy
< \infty \eqno (4.9)$$
Since 
 we are assuming $\eta$ was chosen so that  $\eta - \delta < -\delta_v$, equation $(4.9)$ contradicts Lemma 1.2.  Hence we 
conclude the estimates of parts a) Theorems 1.3 and 1.4 are sharp as desired. 

The proof of sharpness of parts b) of the two theorems is very similar, so we omit the full details. One assumes the result holds for some $\epsilon' > \epsilon$, chooses some $\eta > 0 $ with $\eta + \epsilon
< \epsilon'$ and instead of using $(4.2)$, one uses
$$I_L = \int K_{\phi,x_0,y_0}(t,u) \psi\bigg({|(t,u)| \over L}\bigg)|(t,u)|^{\epsilon ' - 2 - \eta}\,ds \eqno (4.10)$$
Then the steps analogous to $(4.2)-(4.9)$ lead to 
$$\sup_r r^{\eta - \epsilon'} \int_{\{(x,y) \in N:\, r < |(x,y)| < 2r\}} g(x,y)\,dx\,dy
< \infty \eqno (4.11)$$
Since the exponent $\eta - \epsilon '$ is less than $-\epsilon$, equation $(1.5)$ is contradicted and we must have sharpness. This completes the proofs of the sharpness statements of Theorems 1.3 and 1.4.

As for Theorem 1.6, one can readily reduce it to the above sharpness statements. Suppose $m(x,y) = \sum_i m_i(x,y)$ satisfies
the conditions of Theorem 1.6, and $i_0$ is an index such the estimate for $m_{i_0}(x_0 + x,y_0 + y)$ given by Theorem 1.3
or 1.4 is stated to be sharp and such that the estimate for $ K(t,u)$ of $(1.4)$ given by adding the estimates for all
$K_{\phi,x_0,y_0}(t,u)$ over all $i$ is the estimate for the $K_{\phi,x_0,y_0}(t,u)$ corresponding to $i = i_0$. We suppose for
argument's sake the estimate for $m_{i_0}(x_0 + x,y_0 + y)$ derives from part a) of Theorem 1.3; the other cases are dealt with
similarly.

If $K(t,u)$ satisfied a better estimate $|K(t,u)| \leq C(2 + |(t,u)|)^{-\delta}$ on the strip $l_H$ of Theorem 1.3a), where 
$\delta > \delta_v$, then in place of $(4.1)$ we would have $|K(sv)| \leq C(2 + |s|)^{-\delta}$. Instead of defining $I_L$
as in $(4.2)$, for a constant $a$ to be determined, one uses
$$I_L = \int K(s v) \psi\bigg({s \over L}\bigg)e^{-isa}|s|^{\delta - 1 - \eta}\,ds \eqno (4.12)$$
One gets that $\sup_L |I_L| < \infty$ exactly as before. Performing the steps from $(4.2)-(4.9)$ this time leads to
$$\sup_r r^{\eta - \delta} \int_{\{(x,y) \in N:\, r < |((x,y) - av) \cdot v| < 2r\}} m(x,y)\,dx\,dy
< \infty \eqno (4.13)$$
By the assumptions of Theorem 1.6, each $m_i(x, y)$ is nonnegative, so we must also have
$$\sup_r r^{\eta - \delta} \int_{\{(x,y) \in N:\, r < |((x,y) - av) \cdot v| < 2r\}} m_{i_0}(x,y)\,dx\,dy < \infty \eqno (4.14)$$
We choose $a$ so that 
$${\{(x,y) \in N:\, r < |((x + x_0,y + y_0) - av) \cdot v| < 2r\}}$$
$$ = {\{(x,y) \in N:\, r < |(x,y) \cdot v| < 2r\}} \eqno (4.15)$$
 Then changing variables from $(x,y)$ to $(x + x_0,y + y_0)$ in $(4.14)$ leads to 
$$\sup_r r^{\eta - \delta} \int_{\{(x,y) \in N:\, r < |(x,y) \cdot v| < 2r\}} m_{i_0}(x_0 + x, y_0 + y)\,dx\,dy < \infty \eqno (4.16)$$
Since $\alpha(x,y)$ is bounded below by a positive constant, $(4.16)$ implies $(4.9)$, and we get a contradiction like before. Thus the estimate for $K(t,u)$ here is in fact sharp.

Although the above dealt with the situation when the estimate for $m(x,y)$ derives from a sharp estimate for 
$m_{i_0}(x_0 + x,y_0 + y)$ derived from part a) of Theorem 1.3 or 1.4, the other situations are dealt with in the analogous manner.
 This concludes the proof of Theorem 1.6.

\noindent {\heading References.}

\noindent [CDMaM]  M. Cowling, S. Disney, G. Mauceri, and D. Muller {\it Damping oscillatory integrals}, 
Invent. Math. {\bf 101}  (1990),  no. 2, 237-260. \parskip = 4pt\baselineskip = 6pt

\noindent [CMa1] M. Cowling, G. Mauceri, {\it Inequalities for some maximal functions. II},  Trans. Amer. 
Math. Soc. {\bf 298} (1986),  no. 1, 341-365.

\noindent [CMa2] M. Cowling, G. Mauceri, {\it Oscillatory integrals and Fourier transforms of surface carried measures}, 
 Trans. Amer. Math. Soc. 304 (1987), no. 1, 53-68. 

\noindent [G1] M. Greenblatt, {\it Uniform bounds for Fourier transforms of surface measures in $\R^3$ with nonsmooth density}, to appear, Trans. Amer. Math. Soc. 

\noindent [G2] M. Greenblatt, {\it Resolution of singularities in two dimensions and the stability of integrals}, Adv. Math., 
{\bf 226} no. 2 (2011) 1772-1802.

\noindent [G3] M. Greenblatt, {\it Fourier transforms of powers of well-behaved 2D real analytic functions}, preprint.

\noindent [G4] M. Greenblatt, {\it Maximal averages over hypersurfaces and the Newton polyhedron},
J. Funct. Anal. 262 (2012), no. 5, 2314-2348. 

\noindent [G5] M. Greenblatt, {\it Van der Corput lemmas and Fourier transforms of irregular hypersurface measures}, preprint.
arxiv:1409.4059

\noindent [Gr] P. Gressman, {\it Damping oscillatory integrals by the Hessian determinant via Schr\"oding- 
er}, preprint.

\noindent [IM] I. Ikromov, M. Kempe, and D. M\"uller, {\it Damped oscillatory integrals and boundedness of
maximal operators associated to mixed homogeneous hypersurfaces} (English summary) Duke Math. J. {\bf 126} 
(2005), no. 3, 471--490.

\noindent [IoSa1] A. Iosevich, E. Sawyer, {\it Oscillatory integrals and maximal averages over homogene- ous
surfaces}, Duke Math. J. {\bf 82} no. 1 (1996), 103-141.

\noindent [IoSa2] A. Iosevich, E. Sawyer, {\it Maximal averages over surfaces},  Adv. Math. {\bf 132} 
(1997), no. 1, 46--119.

\noindent [PS] D. H. Phong, E. M. Stein, {\it The Newton polyhedron and
oscillatory integral operators}, Acta Mathematica {\bf 179} (1997), 107-152.

\noindent [S] E. Stein, {\it Harmonic analysis; real-variable methods, orthogonality, and oscillatory inte- grals}, Princeton Mathematics Series Vol. 43, Princeton University Press, Princeton, NJ, 1993. 

\noindent [SoS] C. Sogge and E. Stein, {\it Averages of functions over hypersurfaces in $R^n$}, Invent.
Math. {\bf 82} (1985), no. 3, 543--556.\parskip = 12pt \baselineskip = 12pt

\noindent Department of Mathematics, Statistics, and Computer Science \hfill \break
\noindent University of Illinois at Chicago \hfill \break
\noindent 322 Science and Engineering Offices \hfill \break
\noindent 851 S. Morgan Street \hfill \break
\noindent Chicago, IL 60607-7045 \hfill \break
\noindent greenbla@uic.edu \hfill\break

\end